\documentclass[twocolumn]{autart}    % Enable this line and disable the
\usepackage{mathrsfs,amsmath,amssymb,booktabs,array,xcolor,url,cite,color,soul,multirow,multicol,balance}
\usepackage{mathbbold,bbm,bm}
\usepackage{stfloats}
\usepackage[dvips]{graphicx}
\usepackage{textcomp}
\usepackage{ulem}
\usepackage{cases}
\usepackage{flushend,cuted}
\usepackage{amsfonts}
\usepackage{indentfirst}
\newtheorem{theorem}{Theorem}
\newtheorem{corollary}{Corollary}
\newtheorem{definition}{Definition}
\newtheorem{proposition}{Proposition}

\newtheorem{remark}{Remark}

\def\dd{\text{d}}

\makeatletter
\def\@biblabel#1{~}

\leftmargin\labelwidth
\advance\leftmargin\labelsep
\advance\leftmargin by 2em%
\usecounter{enumiv}%

\sethlcolor{yellow}

\makeatother
\begin{document}

\begin{frontmatter}

\title{Stabilization of (state, input)-disturbed CSTRs through the port-Hamiltonian systems approach}%\thanksref{footnoteinfo}} % Title, preferably not more

\thanks[footnoteinfo]{The authors are with the School of Mathematical Sciences, Zhejiang
University, Hangzhou 310027, China (e-mail: 11535029@zju.edu.cn; zhou\_fang@zju.edu.cn; gaochou@zju.edu.cn (Correspondence: C. H. Gao)).}

\author{Yafei~Lu, Zhou~Fang, Chuanhou~Gao}%\ead{cicero@senate.ir},    % Add the
%\author[Rome]{Julius Caesar}\ead{julius@caesar.ir},               % e-mail address
%\author[Baiae]{Publius Maro Vergilius}\ead{vergilius@culture.ir}  % (ead) as shown

%\address[Paestum]{Buckingham Palace, Paestum}  % Please supply
%\address[Rome]{Senate House, Rome}             % full addresses
%\address[Baiae]{The White House, Baiae}        % here.

\begin{keyword}                           % Five to ten keywords,
CSTR, port-Hamiltonian, stochastic passivity, entropy, (state, input)-disturbed% chosen from the IFAC
\end{keyword}                             % keyword list or with the
                                          % help of the Automatica
                                          % keyword wizard

\begin{abstract}
It is a universal phenomenon that the state and input of the continuous stirred tank reactor (CSTR) systems are both disturbed. This paper proposes a (state, input)-disturbed port-Hamiltonian framework that can be used to model and further designs a stochastic passivity based controller to asymptotically stabilize in probability the (state, input)-disturbed CSTR (sidCSTR) systems. The opposite entropy function and the availability function are selected as the Hamiltonian for the model and control purposes, respectively. Furthermore, the proposed (state, input)-disturbed port-Hamiltonian model can simultaneously characterize the first law and the second law of thermodynamics when the opposite entropy function acts as the Hamiltonian. A simple CSTR example illustrates how the port-Hamiltonian method is utilized for modeling and controlling sidCSTR systems.
\end{abstract}
\end{frontmatter}

\section{Introduction}
Continuous stirred tank reactors (CSTRs) are a class of widely used continuous reactors in chemical industrial processes. The operation of a CSTR takes place with reactants fed continuously into the tank through ports at the top while products removed continuously from the bottom, shown in Fig. $1$. A basic assumption for this device is that it admits perfect mixing so that the contents, like temperature, concentrations, etc., have uniform properties throughout the tank, and moreover, the exit stream has the same contents as in the tank. Even so, the CSTR process still exhibits highly complex behaviors due to nonlinear coupling among irreversible thermodynamics, reaction kinetics and hydrodynamics, characterized by multiple equilibrium points, instability, lack of precise physical interpretation, etc. \cite{Favache2009a}. It thus poses a great challenge to investigate CSTRs for both academical research and engineering applications.

\begin{figure}
\begin{center}
{\includegraphics[width=2.3in]{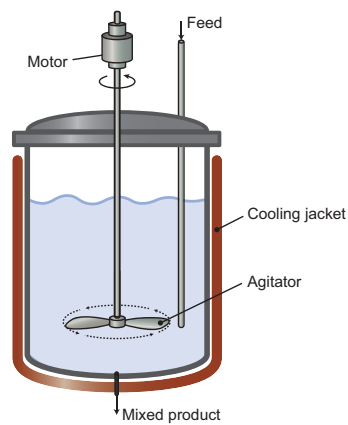}} {\caption{
Schematic diagram of a CSTR.}}
\end{center}
\end{figure}

Sustainable production concepts require the CSTR process to run with optimal performance, which often results from the unstable operation conditions \cite{Hoang2008}. This motivates the stabilization of CSTR around an unstable steady state to become an active research issue during the past decades. The early stabilizing strategies include input/output feedback linearization \cite{Adebekun1991} and nonlinear PI control algorithm \cite{Alvarez-Ramirez1999} while the recent studies focus on linking thermodynamics with mathematical system theory \cite{Ydstie1997, Alonso2001, Ruszkowski2005, Favache2009a}. The purpose of the latter strategy is to develop a entropy-based Lyapunov function, a notion derived from thermodynamics concepts, for transport reaction systems, of which the CSTR is the most typical representative. Ydstie and his coworkers (1997,2001,2005) firstly connected the thermodynamic system with the passivity theory of nonlinear control through taking the Legendre transform of the entropy difference between different states as the storage function, and then developed the inventory (like total mass, energy or the holdup of a chemical constituent) control approach based on the theory of dissipative systems to reach stabilization. Favache and Dochain (2009) performed a detailed study on finding Lyapunov function candidates related to the thermodynamics, including the entropy, the entropy production and the internal energy. Their applicable scopes were also exhibited, such as the entropy only suitable for isolated systems, the entropy production emphasized with linear phenomenological laws, etc. Another class of investigations \cite{Eberard2007, Favache2010, Ramirez2013,Ramírez2016, Hoang2011, Hoang2012} manage to formulate the dynamics of the CSTR into a structure of a port/quasi-port Hamiltonian system with the main motivation of expressing the first and second principle of thermodynamics simultaneously. However, this good imagination sensitively depends on the selection of the Hamiltonian. The negative entropy as the Hamiltonian fails to support the imagination \cite{Hoang2011} while the internal energy based availability function as the Hamiltonian supports it \cite{Ramírez2016}. The stabilization of the CSTR system in the form of the port/quasi-port Hamiltonian structure may be reached through relating the Hamiltonian structure to the passivity control theory \cite{van2000, Alonso2001, van2004}.

Despite extensive studies made towards stabilizing the CSTR systems in the literature, there is no unified method applicable to a wide class of transport reaction systems. Their control as well as the link between thermodynamics and the Lyapunov stability theory remains an open problem \cite{Favache2009a}. The circumstance will get worse when the stochastic phenomena in the CSTR process have to be considered. In fact, it is highly possible that the CSTR process has uncertainties, such as material compositions change, input/output flow fluctuation, external environment temperature fluctuation, reaction kinetics uncertainty, etc. Frankly speaking, the stochastic model can capture the nature of the CSTR process better. However, little attention is paid to stochastic CSTR processes. Even though the kinetic uncertainty is considered in stabilizing the CSTR process, which acts as an example of the stochastic nonlinear time-delay system \cite{Liu2017}, the input disturbance is often ignored. For this reason, the current work is devoted to stabilizing in probability the CSTR systems with state and input both disturbed. By utilizing the (state, input)-disturbed stochastic port-Hamiltonian system (sidSPHS) framework \cite{Fang2017}, we can rewrite the (state, input)-disturbed CSTR (sidCSTR) model to be a sidSPHS with the opposite entropy function as the Hamiltonian. A notable advantage for this rewritten version is that it can characterize both the first law and second law of thermodynamics. Further, by means of the improved stochastic generalized canonical transformation $\mathscr{T}$ \cite{Satoh13,Fang2017}, the rewritten model is transformed to another dynamically equivalent sidSPHS with the availability function as the new Hamiltonian. The latter exhibits stochastic passivity with respect to the new Hamiltonian, based on which a stochastic passivity based controller is developed to asymptotically stabilize the transformed sidSPHS system as well as the original sidCSTR system in probability.

The rest of this paper is organized as follows. Section 2 reviews some basic concepts about stochastic stability and sidSPHS framework. This is followed by the formulation of sidCSTRs into a sidSPHS in Section 3. Section 4 develops a stochastic passivity based controller to asymptotically stabilize sidCSTR systems in probability. A sidCSTR system with first-order reaction acts as a case study to illustrate the proposed stabilizing strategy in Section 5. Finally, Section 6 concludes the paper.

\textit{Notion:} $\mathbb{R}^n$, $\mathbb{R}_\geq0$ and $\mathbb{R}^+$ represent the space of n-dimensional real vectors, the set of nonnegative real numbers and positive real numbers respectively ; $\mathbb{R}^{n\times m}$ represents the set of all ($n \times m$) real matrices ; $\mathbb{Z}_+$ represents the set of positive integers; $\mathbbold{0}_n$ and $\mathbbold{1}_n$ represent n-dimensional vector with all entries equal to 0 and 1, respectively; $\text{E}(\cdot)$ represents mathematical expectation; $|\cdot|$ and $\|\cdot\|$ represent absolute value and 2-norm respectively; $\mathrm{tr}\{{A}\}$, $A_{\cdot i}$ and $A_{j \cdot }$ represent the trace of matrix $A$ , the $i$th column and the $j$th row of matrix $A$ respectively; $\bm{x}^\top$ represents the transpose of the vector (or matrix) $\bm{x}$;  $\mathbbold{1}_{n\times n}$ represents ($n\times n$) matrix with all entries equal to 1; $\mathrm{diag}(x_1, \cdots, x_n)$ represents the ($n\times n$) diagonal matrix whose diagonal elements are $x_j (j=1,\cdots,n)$ while other elements equal to 0; $\nabla (\cdot)$ represents the gradients of $(\cdot)$.

%%%%%%%%%%%%%%%%%%%%%%%%%%%%%%%%%%%%%%%%%%%%%%%%%%%%%%%%%%%%%%%%
\section{(state, input)-disturbed stochastic port-Hamiltonian systems}
{In this section, some basic concepts and results
about stability in the sense of probability and (state, input)-disturbed SPHSs are recalled.}

\subsection{Stability in probability}
Consider the input/output stochastic differential system in $\mathbb{R}^n\times\mathbb{R}^m$ written in the sense of It$\mathrm{\hat{o}}$:
\begin{eqnarray}{\label{SS}}
\left\{
\begin{array}{lll}
\dd\textit{\textbf{x}}=\textit{\textbf{f}}(\textit{\textbf{x}},\textit{\textbf{u}})\dd t+\textit{\textbf{a}}(\textit{\textbf{x}})\dd\boldsymbol{\omega},
\\
\textit{\textbf{y}}=\textit{\textbf{b}}(\textit{\textbf{x}},\textit{\textbf{u}}),
\end{array}
\right.
\end{eqnarray}
where
\begin{enumerate}
\item {$\textit{\textbf{x}}\in {\mathbb{R}^n}$, $\textit{\textbf{u}},~\textit{\textbf{y}}\in\mathbb{R}^m~(m\leq n)$ are the state, measurable input and output, respectively, $t\in\mathbb{R}_{\geq0}$ is the time;}
\item {$\textit{\textbf{f}}:\mathbb{R}^n\times\mathbb{R}^m\mapsto\mathbb{R}^n$, $\textit{\textbf{a}}:\mathbb{R}^n\mapsto \mathbb{R}^{n\times r}$ and $\textit{\textbf{b}}:\mathbb{R}^n\times\mathbb{R}^m\mapsto \mathbb{R}^{m}$ are all local Lipschitz continuous functions such that they vanish in the equilibrium, denoted by $\textit{\textbf{x}}^*$, usually set as the origin, in the unforced case.}
     %and for any $\textit{\textbf{x}} \in {\mathbb{R}^n}$ and $\textit{\textbf{u}} \in {\mathbb{R}^m}$ there is a nonnegative constant $C$ satisfying $$\parallel\textit{\textbf{f}}(\textit{\textbf{x}},\textit{\textbf{u}})\parallel+\parallel \textit{\textbf{a}}(\textit{\textbf{x}})\parallel\leq C(1+\parallel \textit{\textbf{x}}\parallel+\parallel \textit{\textbf{u}}\parallel);$$}
\item {$\boldsymbol{\omega}$ is an $\mathbb{R}^r$-valued standard Wiener process defined on a complete probability space.}
\end{enumerate}

The stability of the equilibrium of the stochastic differential equation (\ref{SS}) in the sense of probability is characterized by the following definition \cite{Khasminskii2011}.

\begin{definition}
The equilibrium solution $\textit{\textbf{x}}^*$ of the stochastic differential equation (\ref{SS}) is
\begin{enumerate}
\item{stable in probability if $\forall \epsilon\textgreater 0$ there is
\begin{equation*}
\lim_{\textit{\textbf{x}}(0)\to
\textit{\textbf{x}}^*} \mathrm{P}(\sup_{t\geq 0}{\parallel\textit{\textbf{x}}(t)-\textit{\textbf{x}}^*\parallel}<\epsilon)=1;
\end{equation*}}
\item{locally asymptotically stable in probability if
\begin{equation*}
\lim_{\textit{\textbf{x}}(0)\to
\textit{\textbf{x}}^*}\mathrm{P}(\lim_{t\to\infty}{\parallel\textit{\textbf{x}}(t)-\textit{\textbf{x}}^*\parallel}=0)=1;
\end{equation*}}
\item{globally asymptotically stable in probability if $\forall
\textit{\textbf{x}}(0)\in\mathbb{R}^n$ there is
\begin{equation*}
\mathrm{P}(\lim_{t\to\infty}{\parallel\textit{\textbf{x}}(t)-\textit{\textbf{x}}^*\parallel}=0)=1,
\end{equation*}}
\end{enumerate}
where $\textit{\textbf{x}}(0)$ is the initial state and $\mathrm{P}(\cdot)$ represents the probability function.
\end{definition}

Besides definition, stochastic Lyapunov theorem is a frequently-used tool to render the stability in probability. Another highly related notation is stochastic passivity \cite{Florchinger1999} that plays an important role in stabilizing stochastic systems, defined below.

\begin{definition}
For an input-output stochastic system in the sense of It$\mathrm{\hat{o}}$, governed by Eq. (\ref{SS}), if there exists a non-negative and twice differentiable storage function $V:\mathbb{R}^n\mapsto\mathbb{R}_{\geq 0}$ such that
\begin{equation}
\mathcal{L}[V(\textbf{x})]\leq \textbf{y}^\top \textbf{u}, \quad\quad\forall\textbf{x}\in\mathbb{R}^n,
\end{equation}
then the system is said to be stochastically passive with respect to $V(\cdot)$ , where $\mathcal{L}[\cdot]$ is the infinitesimal generator of the stochastic process solution $\textit{\textbf{x}}(t)$ of the
stochastic differential equation (\ref{SS}), defined by
\begin{equation*}
\mathcal{L}[V(\textit{\textbf{x}})]=\nabla^\top{V}(\textbf{x})\textit{\textbf{f}}(\textit{\textbf{x}},\textit{\textbf{u}})
+\frac{1}{2}\mathrm{tr}\left\{\frac{\partial^2{V}(\textbf{x})}{\partial \bm{x}^2}\textit{\textbf{a}}(\textit{\textbf{x}})\textit{\textbf{a}}(\textit{\textbf{x}})^\top\right\}.\notag\\
\end{equation*}
\end{definition}

If the storage function $V(\cdot)$ is positive definite with respect to $\textit{\textbf{x}}-\textit{\textbf{x}}^*$, it is easy to realize stabilization of the stochastic system (\ref{SS}) in the sense of probability through connecting it with a proportional controller in negative feedback \cite{Fang2017}.

\subsection{(state, input)-disturbed stochastic port-Hamiltonian systems}
A special case of the stochastic system (\ref{SS}) is the stochastic port-Hamiltonian system proposed by Satoh \& Fujimoto (2013). Essentially, both of them only capture the process noise, but fail to consider the input noise. As an extension, the sidSPHS  \cite{Fang2017} can capture both process noise and input noise, which is of the form
\begin{eqnarray}{\label{sid-SPHS}}
\left\{
\begin{array}{lll}
\dd\textit{\textbf{x}}=\big[(\textit{\textbf{J}}(\textit{\textbf{x}})-\textit{\textbf{R}}(\textit{\textbf{x}}))\nabla H(\textit{\textbf{x}})+\textit{\textbf{g}}(\textit{\textbf{x}})\textit{\textbf{u}}\big]\dd t\\
~~~~~~+\textit{\textbf{a}}(\textit{\textbf{x}})\dd\boldsymbol{\omega}_1+\boldsymbol{\gamma}(\textit{\textbf{x}})\tilde{\textit{\textbf{u}}}\boldsymbol{\sigma}(\textit{\textbf{x}})\dd\boldsymbol{\omega}_2, \\
\textit{\textbf{y}}=\textit{\textbf{g}}^\top(\textit{\textbf{x}})\nabla H(\textit{\textbf{x}})+\boldsymbol{\delta}(\textit{\textbf{x}})\textit{\textbf{u}},
\end{array}
\right.
\end{eqnarray}
where
\begin{enumerate}
\item{$\textit{\textbf{J}},~\textit{\textbf{R}}:\mathbb{R}^n\mapsto\mathbb{R}^{n\times n}$, with $\textit{\textbf{J}}=-\textit{\textbf{J}}^\top,~\textit{\textbf{R}}=\textit{\textbf{R}}^\top\succeq0$ (positive semi-definite), are the interconnection and damping matrices, respectively;}
\item{$\boldsymbol{\omega}_1$ and $\boldsymbol{\omega}_2$ are two mutually independent $\mathbb{R}^{r_{1}}$-valued and $\mathbb{R}^{r_{2}}$-valued standard Wiener processes defined on a probability space, respectively;}
\item{$H:\mathbb{R}^n\mapsto\mathbb{R}$ is the Hamiltonian representing
    the total storied energy, $\textit{\textbf{g}}:\mathbb{R}^n\mapsto\mathbb{R}^{n\times m}$ describes the control port, $\textit{\textbf{u}}=(u_1,\cdots, u_m)^\top$ is the control action with the expectation satisfying $\text{E} (\int_0^t \|\textit{\textbf{u}}(\tau)\|^2 d\tau)\leq \infty$ for all $t\in\mathbb{R}_{\geq 0}$, $\textit{\textbf{a}}:\mathbb{R}^n\mapsto\mathbb{R}^{n\times r_1}$ is the process noise port, and $\boldsymbol{\gamma}(\textit{\textbf{x}})\tilde{\textit{\textbf{u}}}\boldsymbol{\sigma}(\textit{\textbf{x}})$
    with $\boldsymbol{\gamma}:\mathbb{R}^n\mapsto\mathbb{R}^{n\times m},~\tilde{\textit{\textbf{u}}}=\text{diag}(u_1,\cdots,u_m)$ and
    $\boldsymbol{\sigma}:\mathbb{R}^n\mapsto\mathbb{R}^{m\times r_2}$ is the input drift term;}
\item{$\boldsymbol{\delta}(\textit{\textbf{x}})\textit{\textbf{u}}$, with $\boldsymbol{\delta}:\mathbb{R}^n\mapsto\mathbb{R}^{m\times m}$, $\boldsymbol{\delta}\succeq0$ and its Frobenius norm less than $1$, i.e., $\|\boldsymbol{\delta}(\textit{\textbf{x}})\|_\mathcal{F} <1$, measures the contribution of the input to output.}
\end{enumerate}
In addition, the diffusion term $\big[(\textit{\textbf{J}}(\textit{\textbf{x}})-\textit{\textbf{R}}(\textit{\textbf{x}}))\nabla H(\textit{\textbf{x}})+\textit{\textbf{g}}(\textit{\textbf{x}})\textit{\textbf{u}}\big]$ and the drift terms $\textit{\textbf{a}}(\textit{\textbf{x}})$, $\boldsymbol{\gamma}(\textit{\textbf{x}})\tilde{\textit{\textbf{u}}}\boldsymbol{\sigma}(\textit{\textbf{x}})$ also satisfy local Lipschitz continuity and vanish at $\textit{\textbf{x}}^*$ if no driving force is imposed. All these conditions guarantee that the stochastic differential equation (\ref{sid-SPHS}) admits a unique solution given any initial state and the stochastic version of Lyapunov theorem can be applied.

For the structure of sidSPHSs in Eq. (\ref{sid-SPHS}), there are two points that need to be noted. One is that the Hamiltonian $H$ is assumed to be neither positive semi-definite nor bounded from below, the other is that there is an output feedthrough so that this equation looks not belonging to the class of port-Hamiltonian systems. In fact, the form of Eq. (\ref{sid-SPHS}) is indeed a port-Hamiltonian system, which can be verified by interconnecting two sidSPHSs in negative feedback to yield a new sidSPHS.

\begin{proposition}
The negative feedback interconnection of any two $\text{sidSPHSs}$ given by Eq. (\ref{sid-SPHS}) yields another $\text{sidSPHS}$.
\end{proposition}
\noindent{\textbf{Proof.}
We use subscripts ``$1$" and $``2"$ to identify two sidSPHSs, respectively. For simplicity but without loss of generalization, the state, input and output of these two systems are assumed to satisfy $\mathrm{dim}(\textit{\textbf{x}}_1)=\mathrm{dim}(\textit{\textbf{x}}_2)$ and $\mathrm{dim}(\textit{\textbf{u}}_1)=\mathrm{dim}(\textit{\textbf{u}}_2)=\mathrm{dim}(\textit{\textbf{y}}_1)=\mathrm{dim}(\textit{\textbf{y}}_2)$. The interconnection of these two systems in negative feedback follows
\begin{eqnarray*}
\left(
\begin{array}{c}
\textit{\textbf{u}}_{1}\\
\textit{\textbf{u}}_{2}
\end{array}
\right)
&=&
\left(
\begin{array}{cc}
\mathbbold{0} & -\textit{\textbf{I}}\\
\textit{\textbf{I}} & \mathbbold{0}
\end{array}
\right)
\left(
\begin{array}{c}
\textit{\textbf{y}}_{1}\\
\textit{\textbf{y}}_{2}
\end{array}
\right)\notag \\
&=&
\left(
\begin{array}{cc}
\mathbbold{0} & -\textit{\textbf{I}}\\
\textit{\textbf{I}} & \mathbbold{0}
\end{array}
\right)
\left(
\begin{array}{c}
\textit{\textbf{g}}_{1}(\textit{\textbf{x}}_1)^\top\frac{\partial
H_{1}(\textit{\textbf{x}}_1)}{\partial
\textit{\textbf{x}}_{1}}+\boldsymbol{\delta}_{1}(\textit{\textbf{x}}_1)\textit{\textbf{u}}_{1}\\
\textit{\textbf{g}}_{2}(\textit{\textbf{x}}_2)^\top\frac{\partial
H_{2}(\textit{\textbf{x}}_2)}{\partial
\textit{\textbf{x}}_{2}}+\boldsymbol{\delta}_{2}(\textit{\textbf{x}}_2)\textit{\textbf{u}}_{2}
\end{array}
\right),
\end{eqnarray*}
where $\textit{\textbf{I}}$ is the identity matrix with suitable dimension. Since
$\Vert\boldsymbol{\delta}_1(\textit{\textbf{x}}_1)\Vert_{\mathcal{F}}\textless 1$ and $\Vert\boldsymbol{\delta}_2(\textit{\textbf{x}}_2)\Vert_{\mathcal{F}}\textless 1$, we have $$\Vert\boldsymbol{\delta_1}(\textit{\textbf{x}}_1)\boldsymbol{\delta}_2(\textit{\textbf{x}}_2)\Vert_{\mathcal{F}}\leq \Vert\boldsymbol{\delta}_1(\textit{\textbf{x}}_1)\Vert_{\mathcal{F}}\Vert\boldsymbol{\delta}_2(\textit{\textbf{x}}_2)\Vert_{\mathcal{F}}\textless 1.$$ Hence, $\textit{\textbf{I}}+\boldsymbol{\delta_1}(\textit{\textbf{x}}_1)\boldsymbol{\delta}_2(\textit{\textbf{x}}_2)$ and $\textit{\textbf{I}}+\boldsymbol{\delta}_2(\textit{\textbf{x}}_2)\boldsymbol{\delta}_1(\textit{\textbf{x}}_1)$ are both nonsingular. The interconnection equation changes to be
\begin{eqnarray*}{\label{input}}
\left(
\begin{array}{c}
\textit{\textbf{u}}_{1}\\
\textit{\textbf{u}}_{2}
\end{array}
\right)
&=& \left(
\begin{array}{cc}
-\boldsymbol{\delta}_{2}\big(\textit{\textbf{I}}+\boldsymbol{\delta}_{1}\boldsymbol{\delta}_{2}\big)^{-1} & -\big(\textit{\textbf{I}}+\boldsymbol{\delta}_{2}\boldsymbol{\delta}_{1}\big)^{-1}\\
\big(\textit{\textbf{I}}+\boldsymbol{\delta}_{1}\boldsymbol{\delta}_{2}\big)^{-1}
&
-\boldsymbol{\delta}_{1}\big(\textit{\textbf{I}}+\boldsymbol{\delta}_{2}\boldsymbol{\delta}_{1}\big)^{-1}
\end{array}
\right) \notag\\
&\times& \left(
\begin{array}{cc}
\textit{\textbf{g}}_{1}(\textit{\textbf{x}}_1)^\top&\mathbbold{0}\\
 \mathbbold{0}& \textit{\textbf{g}}_{2}(\textit{\textbf{x}}_2)^\top
\end{array}
\right)
\left(
\begin{array}{c}
\frac{\partial
H_{1}(\textit{\textbf{x}}_1)}{\partial
\textit{\textbf{x}}_{1}}\\
\frac{\partial
H_{2}(\textit{\textbf{x}}_2)}{\partial
\textit{\textbf{x}}_{2}}
\end{array}
\right).
\end{eqnarray*}
Combining this equation with Eq.(\ref{sid-SPHS}) yields
\begin{eqnarray*}
\dd
\left(
\begin{array}{c}
\textit{\textbf{x}}_{1}\\
\textit{\textbf{x}}_{2}
\end{array}
\right) &=& \left[ \tilde{\bm{J}}- \tilde{\bm{R}}\right]
\left(
\begin{array}{c}
\frac{\partial
H_{1}}{\partial
\textit{\textbf{x}}_{1}}\\
\frac{\partial
H_{2}}{\partial
\textit{\textbf{x}}_{2}}
\end{array}
\right)\dd t+
\left(
\begin{array}{c}
\textit{\textbf{a}}(\textit{\textbf{x}}_1)\\
\textit{\textbf{a}}(\textit{\textbf{x}}_2)
\end{array}
\right)
\dd\boldsymbol{\omega}_1 \notag\\
&+&
\left(
\begin{array}{cc}
\boldsymbol{\gamma}(\textit{\textbf{x}}_1)&0\\
0&\boldsymbol{\gamma}(\textit{\textbf{x}}_2)
\end{array}
\right)
\left(
\begin{array}{cc}
\tilde{\bm{u}}_1&0\\
0&\tilde{\bm{u}}_2
\end{array}
\right)
\left(
\begin{array}{c}
\boldsymbol{\sigma}(\textit{\textbf{x}}_1)\\
\boldsymbol{\sigma}(\textit{\textbf{x}}_2)
\end{array}
\right)
\dd\boldsymbol{\omega}_2
\end{eqnarray*}
with
\begin{eqnarray*}
\tilde{\bm{J}}&= &
\left(
\begin{array}{cc}
 \textit{\textbf{g}}_{1}(\textit{\textbf{x}}_1)&\mathbbold{0}\\
 \mathbbold{0}& \textit{\textbf{g}}_{2}(\textit{\textbf{x}}_2)
\end{array}
\right)
\left(
\begin{array}{cc}
 \mathbbold{0}& -\big(\textit{\textbf{I}}+\boldsymbol{\delta}_{2}\boldsymbol{\delta}_{1}\big)^{-1} \\
\big(\textit{\textbf{I}}+\boldsymbol{\delta}_{1}\boldsymbol{\delta}_{2}\big)^{-1}
&\mathbbold{0}
\end{array}
\right)\notag\\
&\times&
\left(
\begin{array}{cc}
 \textit{\textbf{g}}_{1}^\top(\textit{\textbf{x}}_1)&\mathbbold{0}\\
 \mathbbold{0}&\textit{\textbf{g}}_{2}^\top(\textit{\textbf{x}}_2)
\end{array}
\right)
+\left(
\begin{array}{cc}
\textit{\textbf{J}}_{1}(\textit{\textbf{x}}_1) &\mathbbold{0} \\
 \mathbbold{0}& \textit{\textbf{J}}_{2}(\textit{\textbf{x}}_2)
\end{array}
\right)
\end{eqnarray*}
and
\begin{eqnarray*}
\tilde{\bm{R}} &=&
\left(
\begin{array}{cc}
 \textit{\textbf{g}}_{1}(\textit{\textbf{x}}_1)& \mathbbold{0}\\
 \mathbbold{0}& \textit{\textbf{g}}_{2}(\textit{\textbf{x}}_2)
\end{array}
\right)
\left(
\begin{array}{cc}
\boldsymbol{\delta}_{2}\big(\textit{\textbf{I}}+\boldsymbol{\delta}_{1}\boldsymbol{\delta}_{2}\big)^{-1}  & \mathbbold{0}\\
\mathbbold{0}&\boldsymbol{\delta}_{1}\big(\textit{\textbf{I}}+\boldsymbol{\delta}_{2}\boldsymbol{\delta}_{1}\big)^{-1}
\end{array}
\right)\notag\\
&\times&
\left(
\begin{array}{cc}
 \textit{\textbf{g}}_{1}^\top(\textit{\textbf{x}}_1)&\mathbbold{0}\\
\mathbbold{0}& \textit{\textbf{g}}_{2}^\top(\textit{\textbf{x}}_2)
\end{array}
\right)
+\left(
\begin{array}{cc}
\textit{\textbf{R}}_{1}(\textit{\textbf{x}}_1) &\mathbbold{0}\\
\mathbbold{0}& \textit{\textbf{R}}_{2}(\textit{\textbf{x}}_2)
\end{array}
\right).
\end{eqnarray*}
Note that $\boldsymbol{\delta}_1(\textit{\textbf{x}}_1)=\boldsymbol{\delta}_1^\top(\textit{\textbf{x}}_1)$ and
$\boldsymbol{\delta}_2(\textit{\textbf{x}}_2)=\boldsymbol{\delta}_2^\top(\textit{\textbf{x}}_2)$, we get $\tilde{\bm{J}}=-\tilde{\bm{J}}^\top$,
$\tilde{\bm{R}}=\tilde{\bm{R}}^\top\succeq 0$.
Namely, the above interconnection yields another $\text{sidSPHS}$ with the Hamiltonian to
be $H_1(\textit{\textbf{x}}_1)+H_2(\textit{\textbf{x}}_2)$.
$\Box$

As for the first point related to the Hamiltonian, if it is bounded from below, the sidSPHS in the form of Eq. (\ref{sid-SPHS}) can reach power balance, or be stochastically passive with respect to the Hamiltonian with some moderate conditions added.

\begin{theorem}\cite{Fang2017}
A $\text{sidSPHS}$ described by Eq. (\ref{sid-SPHS}) is stochastically passive with respect to $H(\textit{\textbf{x}})$ if and only if $\forall
\textbf{x}$ the Hamiltonian satisfies $H(\textbf{x})\geq
H(\mathbbold{0}_n)=0$,
\begin{equation}{\label{passivityCond0}}
\frac{1}{2}\mathrm{tr}\left\{\frac{\partial^{2} H(\textit{\textbf{x}})}{\partial
\textbf{x}^{2}}\textbf{a}(\textit{\textbf{x}})\textbf{a}^\top(\textit{\textbf{x}})\right\} \leq
\frac{\partial H(\textit{\textbf{x}})}{\partial
\textbf{x}}^\top\textbf{R}(\textit{\textbf{x}})\frac{\partial H(\textit{\textbf{x}})}{\partial
\textbf{x}}
\end{equation}
and
\begin{equation}\label{sid-SPHSPassive}
\boldsymbol{\delta}(\textbf{x})-\frac{1}{2}\boldsymbol{\sigma}(\textit{\textbf{x}})\boldsymbol{\sigma}^\top(\textit{\textbf{x}})
\circ\boldsymbol{\gamma}^\top(\textit{\textbf{x}})\frac{\partial^{2}{H}(\textbf{x})}{\partial\textbf{x}^{2}}\boldsymbol{\gamma}(\textit{\textbf{x}})\succeq 0,
\end{equation}
where $``\circ"$ is the Hadamard product.
\end{theorem}

Based on \textbf{Theorem 1}, if the Hamiltonian is further assumed to be positive definite, then a simple unity feedback controller $\textit{\textbf{u}}=-\textit{\textbf{y}}$
is able to stabilize the system in probability at the origin. However, the Hamiltonian usually cannot behave as a Lyapunov function, and even not bounded from below. To address this issue, Fang \& Gao (2017), based on the stochastic generalized canonical transformation \cite{Satoh13}, proposed the following transformation
\begin{eqnarray}{\label{ImSGCT}}
\mathscr{T}:~\left\{
\begin{array}{lll}
\bar{\bm{x}}=\varphi(\textit{\textbf{x}}),\\
\bar{H}(\bar{\bm{x}})=\left.H(\textit{\textbf{x}})+H'(\textit{\textbf{x}})\right|_{\textit{\textbf{x}}=\varphi^{-1}(\bar{\bm{x}})},\\
\bar{\bm{y}}_d=\left.\textit{\textbf{y}}_d+\boldsymbol{\alpha}(\textit{\textbf{x}})\right|_{\textit{\textbf{x}}=\varphi^{-1}(\bar{\bm{x}})},\\
\bar{\bm{u}}=\left.\textit{\textbf{u}}+\boldsymbol{\beta}(\textit{\textbf{x}})\right|_{\textit{\textbf{x}}=\varphi^{-1}(\bar{\bm{x}})}
\end{array}
\right.
\end{eqnarray}
to change the Hamiltonian, where $\varphi:\mathbb{R}^n\mapsto \mathbb{R}^n$,
$H':\mathbb{R}^n\mapsto \mathbb{R}$ and
$\boldsymbol{\alpha},~\boldsymbol{\beta}:\mathbb{R}^n\mapsto
\mathbb{R}^m$ are all differentiable functions,  $\boldsymbol{\alpha}(\textit{\textbf{x}})=\textit{\textbf{g}}^\top(\textit{\textbf{x}})\frac{\partial {H'}(\textit{\textbf{x}})}{\partial \textit{\textbf{x}}}$, $\textit{\textbf{y}}_d=\textit{\textbf{g}}^\top(\textit{\textbf{x}})\frac{\partial {
H}(\textit{\textbf{x}})}{\partial \textit{\textbf{x}}}$ and moreover,
$\varphi(\cdot)$ is a diffeomorphism and can preserve the
original point. A lower bounded and even positive definite Hamiltonian may be generated through $\mathscr{T}$.

%%%%%%%%%%%%%%%%%%%%%%%%%%%%%%%%%%%%%%%%%%%%%%%%%%%%%%%%%%%%%%%%%%%%%%%%%%%%%%%
\section{Stochastic Hamiltonian formulation of (state, input)-disturbed CSTRs}
In this section, we use the sidSPHS framework to model CSTRs with state and input disturbed.

Consider a CSTR involving in $p$ chemical components
$X_1,\cdots,X_p$ and $l$ reactions
$\mathscr{R}_1,\cdots,\mathscr{R}_l$ with the $i$th reaction written
as
\begin{equation*}
\mathscr{R}_i:~~~\sum^{p}_{j=1}z_{ji}X_j \leftrightarrows
\sum^{p}_{j=1}z'_{ji}X_j,
\end{equation*}
where $z_{ji}, ~z'_{ji}\in\mathbb{Z}_{\geq 0}$ are the stoichiometric
coefficients. To control the temperature $T$, the CSTR is equipped
with a surrounding jacket, in which the cooling/heating fluid flows
for heat exchange. The dynamics can be deduced from the mass and
energy balances by considering the following assumptions:

\textit{A.1:} The reaction mixture is ideal and incompressible; and
moreover, the reaction volume $V$ and pressure $P$ are supposed to
be constants.

\textit{A.2:} At the inlet the reactor, the pure components $X_j$
are fed at temperature $T_j^{\text{in}}$ and molar concentration
$c_j^{\text{in}}$.

\textit{A.3:} The reaction kinetics obey the Arrhenius law, with
which the forward and backward reaction rate of $\mathscr{R}_i$ are
given by
\begin{eqnarray*}\label{Kinetics}
~\left\{
\begin{array}{lll}
r_{f_i}&=k_{0f_i}\text{exp}\left(-\frac{E_{f_i}}{RT}\right)\Pi_{j=1}^p
c_j^{z_{ji}}, \\
r_{b_i}&=k_{0b_i}\text{exp}\left(-\frac{E_{b_i}}{RT}\right)\Pi_{j=1}^p
c_j^{z'_{ji}},
\end{array}
\right.
\end{eqnarray*}
where $E_{f_i},E_{b_i}$ are reaction activation energy, $k_{0f_i},
k_{0b_i}$ are reaction rate constants, and $c_j$ is the molar
concentration of component $X_j$, evaluated by the molar mass $N_j$
according to $c_j=\frac{N_j}{V}$.

\textit{A.4:} The heat exchange with the jacket is proportional to
the difference between the jacket temperature $T_w$ and the mixture
temperature $T$ in the reactor, i.e, $\dot Q=\lambda (T_w-T)$ with
$\lambda$ the heat transfer coefficient.

\textit{A.5:} There are slight fluctuations on chemical
reactions, inlet/outlet volume flows and heat exchange. All fluctuations are supposed to be mutually independent standard Wiener processes, labeled by $\omega_1$, $\omega_2$,
$\omega_3\in\mathbb{R}$, respectively. They act on the corresponding
objects in proportion to their respective standard deviations $\rho_1, ~\rho_2, ~\rho_3 \in\mathbb{R}^+$, which are assumed to be sufficiently small in comparison to the above mentioned macroscopic variables.
%%%%%%%%%%%%%%%%%%%%%%%%%%%%%%%%%%%%%%%%%%%%%%%%%%%%%%%%%%%%%%%%%%%%%%%%%%%%%%%%%%%%%%%%%%%%%%%%%%%%%%%%%%%%%%%%%%%%%%%%%%%%%%%%%%%%%%%%%%%%%%%%
\newcounter{TempEqCnt}                         % 创建临时变量TempEqCnt
\setcounter{equation}{15}
 % 当前公式序号变为x，x等于长公式应有的序号减1.
\begin{figure*}[!b]
\hrulefill
\begin{equation}{\label {sid-CSTR}}
\begin{split}
~{\Sigma_s:~}\left\{
\begin{array}{lll}
\dd{\textit{\textbf{x}}}&=& (\textit{\textbf{J}}(\bm{x})-\textit{\textbf{R}}(\bm{x}))
\frac{\partial (-S)}{\partial\textit{\textbf{x}}}\dd{t}
+\underbrace{\left(
\begin{array}{cc}
    \frac{U^{\text{in}}}{V}-\frac{U}{V}& 1 \\
     \textit{\textbf{c}}^{\text{in}}-\frac{\textit{\textbf{N}}}{V}& 0 \\
\end{array}\right)}_{{\textit{\textbf{g}}}(\textit{\textbf{x}})}
\underbrace{\left(
\begin{array}{cc}
  q \\
    \dot Q \\
\end{array}\right)}_{{\textit{\textbf{u}}}} \dd{t}
+\underbrace{\left(
\begin{array}{c}
0\\
   (\textit{\textbf{z}}-\textit{\textbf{z}}')(\textit{\textbf{r}}_b-\textit{\textbf{r}}_f)\rho_1 V\\
  \end{array}\right)}_{\bm{a}(\bm{x})}\dd{\omega_1}\\
&+&\underbrace{\left(
\begin{array}{cc}
     \frac{U^{\text{in}}}{V}-\frac{U}{V}& 1 \\
     \textit{\textbf{c}}^{\text{in}}-\frac{\textit{\textbf{N}}}{V}& 0 \\
  \end{array}\right)}_{{\boldsymbol{\gamma}(\bm{x})}}
  \underbrace{  \left(
   \begin{array}{cc}
  q  & 0\\
   0& \dot Q \\
\end{array}\right)}_{{\tilde{\textit{\textbf{u}}}}}
\underbrace{ \left(
\begin{array}{cc}
\rho_2& 0\\
 0 & \rho_3\\
  \end{array}
\right)}_{{\boldsymbol{\sigma}(\textit{\textbf{x}})}}
\left(
\begin{array}{c}
\dd{\omega_2}\\
\dd{\omega_3}\\
  \end{array}
\right)\\
\textit{\textbf{y}}&=&\underbrace{\left(
\begin{array}{ccc}
 \frac{U^{\text{in}}}{V}-\frac{U}{V} & (\textit{\textbf{c}}^{\text{in}}-\frac{\textit{\textbf{N}}}{V})^\top \\
 1&0
\end{array}
\right)}_{{\textit{\textbf{g}}}^\top(\textit{\textbf{x}})}\frac{\partial (-S)}{\partial\textit{\textbf{x}}}
+
\underbrace{\left(
\begin{array}{cc}
\frac{1}{2}\rho^{2}_2 \frac{M}{\theta}&0\\
0&\frac{1}{2}\rho^{2}_3\frac{1}{\theta}
\end{array}
\right)}_
{\boldsymbol{\delta}(\bm{x})} \underbrace{\left(
 \begin{array}{c}
 q  \\
    \dot Q \\
\end{array}\right)}_{{\textit{\textbf{u}}}}
\end{array}
\right.
\end{split}
\end{equation}
 \hrulefill
 \end{figure*}
\setcounter{equation}{6}
%%%%%%%%%%%%%%%%%%%%%%%%%%%%%%%%%%%%%%%%%%%%%%%%%%%%%%%%%%%%%%%%%%%%%%%%%%%%%%%%%%%%%%%%%%%%%%%%%%%%%%%%%%%%%%%%%%%%%%%%%%%%%%%%%%%%%%%%%%%%%%%%%

The mass balances in the It\^{o} form are then given by
\begin{eqnarray}\label{MassBalance}
\dd{\textit{\textbf{N}}}&=&(\textit{\textbf{z}}-\textit{\textbf{z}}')(\textit{\textbf{r}}_b-\textit{\textbf{r}}_f)V(\dd t+\rho_1 \dd\omega_1)\notag\\
&+&\left(q\textit{\textbf{c}}^{in}-q\frac{\textit{\textbf{N}}}{V}\right)(\dd t+\rho_2\dd\omega_2),
\end{eqnarray}
where $\textit{\textbf{N}}=(N_1,\cdots,N_p)^\top$,
$\textit{\textbf{c}}^{\text{in}}=(c_1^{\text{in}},\cdots,c_p^{\text{in}})^\top$,
$\textit{\textbf{z}}=[z_{ji}]_{p\times l}$, $\textit{\textbf{z}}'=[z'_{ji}]_{p\times
l}$, $\textit{\textbf{r}}_f=(r_{f_1},\cdots,r_{f_l})^\top$ and $\textit{\textbf{r}}_b=(r_{b_1},\cdots,r_{b_l})^\top$.
The energy balance, evaluated by the internal energy $U$, in the sense of It\^{o} is of the form
\begin{eqnarray}\label{EnergyBalance}
\dd U&=&\left(q\frac{U^{\text{in}}}{V}-q\frac{U}{V}\right)(\dd t+\rho_2\dd\omega_2)+\dot
Q(\dd t+\rho_3\dd\omega_3),
\end{eqnarray}
where $U^{\text{in}}$ is the inlet internal energy.

A key point to formulate the $\text{sidCSTR}$ system
into a $\text{sidSPHS}$ structure is to select the Hamiltonian. Similar to
the Hamiltonian formulation of the deterministic CSTRs
\cite{Hoang2011}, the opposite entropy $-S$ is chosen as the
Hamiltonian, which follows the famous Gibbs' equation
\begin{equation}\label{Gibbs}
\dd{S}=\frac{1}{T}\dd{U}+\frac{P}{T}\dd{V}-\frac{\boldsymbol{\mu}^\top}{T}\dd{\textit{\textbf{N}}},
\end{equation}
where $\boldsymbol{\mu}=(\mu_1,\cdots,\mu_p)^\top$ is the chemical
potential vector. The entropy function is set to abide by some basic assumptions \cite{Herbert1985}: i) principle of maximum entropy; ii) positive homogeneousness of degree one. From these two assumptions, it is easy to obtain that the entropy is concave. We further assume that the entropy function is twice continuously differentiable, so the Hessian matrix of $-S$ is positive semi-definite. For the sidCSTRs under consideration, note that $\dd V=0$ in Eq. (\ref{Gibbs}), then by letting
$\textit{\textbf{x}}=(U,\boldsymbol{\textit{\textbf{N}}}^\top)^\top$ we have
\begin{equation}\label{EntropyGradient}
\frac{\partial (-S)}{\partial
\textit{\textbf{x}}}=\left(-\frac{1}{T},\frac{\boldsymbol{\mu}^\top}{T}\right)^\top.
\end{equation}
%where $U^*$ is the equilibrium internal energy and $\boldsymbol{\textit{\textbf{N}}}^*$ is the equilibrium molar mass vector.
The calculation of the second derivative of the negative entropy requires the definitions of the heat capacity and chemical potential in the case of the ideal and incompressible mixture \cite{Hoang2011}, i.e.,
\begin{eqnarray}\label{ChemicalPotential}
~\left\{
\begin{array}{lll}
U&=&\bm{N}^\top[\bm{C}_P(T-T_{\mathrm{ref}})+\bm{h}_{\mathrm{ref}}]-PV,\\
\frac{\boldsymbol{\mu}}{T}&=&-\bm{C}_P\ln\frac{T}{T_{\mathrm{ref}}}+R\mathrm{Ln}\left(\frac{\textit{\textbf{N}}}{\textit{\textbf{N}}^\top\mathbbold{1}_{p}}\right)-\bm{s}_{\mathrm{ref}}+\frac{\bm{h}}{T},
\end{array}
\right.
\end{eqnarray}
where $\bm{C}_P=(C_{P_1},\cdots, C_{P_p})^\top$ is the isobaric heat capacity vector, $T_{\mathrm{ref}}$ the reference temperature, $\bm{h}_{\mathrm{ref}}=(h_{1_{\mathrm{ref}}},\cdots,h_{p_{\mathrm{ref}}})^\top$ the reference molar enthalpy vector, $\bm{s}_{\mathrm{ref}}=(s_{1_{\mathrm{ref}}},\cdots,s_{p_{\mathrm{ref}}})^\top$ the reference molar entropy vector and $\bm{h}=\bm{C}_P(T-T_{\mathrm{ref}})+\bm{h}_{\mathrm{ref}}$ represents the molar enthalpy vector. Combining Eqs. (\ref{EntropyGradient}) and (\ref{ChemicalPotential}), we get
\begin{eqnarray}\label{twice entropy}
\frac{\partial^2(-S)}{\partial{\bm{x}^2}}&=&\left(
\begin{array}{cc}
\frac{1}{\theta},&-\frac{\bm{h}^\top}{\theta}\\
-\frac{\bm{h}}{\theta },&\frac{\bm{h}\bm{h}^\top}{\theta}-\frac{R}{\textit{\textbf{N}}^\top\mathbbold{1}_{p}}\mathbbold{1}_{p\times p}+\mathrm{diag}\left(\frac{R}{N_j}\right)
\end{array}
\right)
\end{eqnarray}
with $\theta=T^2\bm{N}^\top\bm{C}_P$, $j=1,\cdots,p$. Throughout the paper, we sometimes drop $``(\textit{\textbf{x}})"$ for notational simplicity.

The following task is to define $\textit{\textbf{J}}$, $\textit{\textbf{R}}$, $\textit{\textbf{u}}$ and others so that Eqs. (\ref{MassBalance}), (\ref{EnergyBalance}) are rewritten according to the sidSPHS structure of Eq. (\ref{sid-SPHS}).

\begin{proposition}
For any sidCSTR governed by Eqs.
(\ref{MassBalance}) and (\ref{EnergyBalance}), assume that the disturbances originating from the inlet/outlet volume flows and heat exchange are weak enough so that
\begin{equation}\label{norm cond}
\rho^{4}_{2}M^2+\rho^{4}_{3}<4\theta^2,
\end{equation}
where
\begin{eqnarray}
M&=&\Delta^\top \bm{c} \left(\bm{h}\bm{h}^\top-\frac{\theta R}{\textit{\textbf{N}}^\top\mathbbold{1}_{p}}\mathbbold{1}_{p\times p}+\mathrm{diag}\left(\frac{\theta R}{N_j}\right)\right)\Delta \bm{c}\notag\\
&-&2\frac{U^{\mathrm{in}}-U}{V}\bm{h}^\top\Delta \bm{c}
+\left(\frac{U^{\mathrm{in}}-U}{V}\right)^2,
\end{eqnarray}
$j=1,\cdots,p$ and $\Delta \bm{c}=\bm{c}^{\mathrm{in}}-\frac{\bm{N}}{V}$. Then the sidCSTR is a sidSPHS, written by Eq. (12) with $\textbf{J}=\mathbbold{0}_{(p+1)\times(p+1)}$, $H=-S$ and
\begin{equation}\label{RMatrix}
\textbf{R}=\left(
\begin{array}{cc}
0,&\mathbbold{0}_{1\times p}\\
\mathbbold{0}_{p\times 1},&VT\sum^{l}_{i=1}\frac{(r_{f_i}-r_{b_i})}{\Delta^\top\textbf{z}_{\cdot i} \boldsymbol{\mu}}
\Delta\textbf{z}_{\cdot i}\Delta^\top\textbf{z}_{\cdot i}
\end{array}
\right),
\end{equation}
where $\Delta\textbf{z}_{\cdot i}=\textbf{z}_{\cdot i}-\textbf{z}'_{\cdot i}$.
\setcounter{equation}{6}
\end{proposition}

\noindent{\textbf{Proof.}}
The proof may be implemented in two separate procedures. The first one is to verify Eq. (\ref{sid-CSTR}) is equivalent to Eqs. (\ref{MassBalance}) and (\ref{EnergyBalance}) while the second one is to prove that the variables defined in Eq. (\ref{sid-CSTR}) satisfies the constraints made in Eq. (\ref{sid-SPHS}).

(i) By inserting Eqs. (\ref{EntropyGradient}) and (\ref{RMatrix}) into Eq. (\ref{sid-CSTR}), it is easy to verify that Eq. (\ref{sid-CSTR}) is completely equivalent to Eqs. (\ref{MassBalance}) and (\ref{EnergyBalance}).

(ii) The sidSPHS structure requests $\bm{R}(\bm{x})$, $\boldsymbol{\delta}(\bm{x})$ to be positive semi-definite and $\|\boldsymbol{\delta}(\bm{x})\|_{\mathcal{F}}\textless 1$.

Let $\boldsymbol{\lambda}=(\lambda_0,\boldsymbol{\lambda}_p^\top)^\top$ be any $\mathbb{R}^{p+1}$ vector, then we have
\begin{eqnarray*}
\quad\quad\boldsymbol{\lambda}^\top\bm{R}\boldsymbol{\lambda}&=&\sum^l_{i=1}\frac{(r_{f_i}-r_{b_i})VT}{\Delta^\top\bm{z}_{\cdot i} \boldsymbol{\mu}} \boldsymbol{\lambda}_p^\top \Delta\bm{z}_{\cdot i}\Delta^\top\bm{z}_{\cdot i}\boldsymbol{\lambda}_p\\
&=&\sum^l_{i=1}\frac{(r_{f_i}-r_{b_i})VT}{\Delta^\top\bm{z}_{\cdot i} \boldsymbol{\mu}}(\Delta^\top\bm{z}_{\cdot i}\boldsymbol{\lambda}_p)^2.
\end{eqnarray*}
Note that $\frac{(r_{f_i}-r_{b_i})VT}{\Delta^\top\bm{z}_{\cdot i} \boldsymbol{\mu}}$ is a well defined nonnegative function since it has the same sign as the entropy creation, $\frac{(r_{f_i}-r_{b_i})V}{T}\Delta^\top\bm{z}_{\cdot i} \boldsymbol{\mu}\geq 0$, of the $i$th chemical reaction \cite{Prigogine1954}, so $\boldsymbol{\lambda}^\top\bm{R}\boldsymbol{\lambda}\geq 0$. Clearly, $\bm{R}=\bm{R}^\top$. We thus get that $\bm{R}$ is positive semi-definite.

The positive semi-definiteness of $\boldsymbol{\delta}(\bm{x})$ may be easily verified from the facts that
$\theta\textgreater 0$, $M=\theta\boldsymbol{\gamma}_{\cdot 1}^\top\frac{\partial^{2}(-S)}{\partial\bm{x}^{2}}\boldsymbol{\gamma}_{\cdot 1}\geq 0$ and $\boldsymbol{\delta}(\bm{x})=\boldsymbol{\delta}^\top(\bm{x})$. Also, since $$\|\boldsymbol{\delta}(\bm{x})\|_\mathcal{F}=\frac{1}{4}\rho_2^4\frac{M^2}{\theta^2}+\frac{1}{4}\rho_3^4\frac{1}{\theta^2},$$ which combines with Eq. (\ref{norm cond}), we obtain $\|\boldsymbol{\delta}(\bm{x})\|_\mathcal{F}\textless 1$. In addition, the noise ports in Eq. (\ref{sid-CSTR}) will vanish in the equilibrium with zero input. $\Box$

\begin{remark}
Under the conditions of no disturbances, the sidCSTR of Eq. (\ref{sid-CSTR}) will degenerate to the deterministic CSTR. If no environmental exchanges occur again, the deterministic model simplifies to be $$\dot{\textit{\textbf{x}}}= (\textit{\textbf{J}}(\bm{x})-\textit{\textbf{R}}(\bm{x}))
\frac{\partial (-S)}{\partial\textit{\textbf{x}}}.$$
This structure can represent the first law and the second law of thermodynamics simultaneously since
\begin{eqnarray*}
~\left\{
\begin{array}{lll}
\frac{\dd U}{\dd t}&=&0,\\
\frac{\dd S}{\dd t}&=&\left(\frac{\partial S}{\partial \textit{\textbf{x}}}\right)^\top\textit{\textbf{R}}(\bm{x})\frac{\partial S}{\partial \textit{\textbf{x}}}\geq 0.
\end{array}
\right.
\end{eqnarray*}
\end{remark}

\begin{remark}
For a sidCSTR, whether the condition of Eq. (\ref{norm cond}) is true mainly depends on the first term $\rho^{4}_{2}M^2$ since the second one $\rho^{4}_{3}$ is far less than $4\theta^2$ according to the assumption A.5. Basically, as long as $\rho^{2}_{2}\boldsymbol{\gamma}_{\cdot 1}^\top\frac{\partial^{2}(-S)}{\partial\bm{x}^{2}}\boldsymbol{\gamma}_{\cdot 1}\textless 2$, Eq. (\ref{norm cond}) will hold almost surely. This may be easily achieved in practice if the temperature and the concentrations of input species are manipulated properly.
\end{remark}

%%%%%%%%%%%%%%%%%%%%%%%%%%%%%%%%%%%%%%%%%%%%%%%%%%%%%%%%%%%%%%%%%%%%%%%%%%%%%%%%%%%%%%%%%%%%%%%%%%%%%%%%%%%%%%%%%%%%%%%%%%%%%%%%%%%%%%%%%%%%%%%%%%
\newcounter{TempEqCnt1}                         % 创建临时变量TempEqCnt
\setcounter{equation}{17}% 当前公式序号变为x，x等于长公式应有的序号减1.
\begin{figure*}[!b]
\hrulefill
\begin{equation}{\label{sid-CSTR1}}
\begin{split}
~{\mathscr{T}(\Sigma_s):~}\left\{
\begin{array}{lll}
\dd{\bar{\bm{x}}}&= &(\bar{\bm{J}}(\bar{\bm{x}})-\bar{\bm{R}}(\bar{\bm{x}}))
\frac{\partial\bar{H}(\bar{\bm{x}})}{\partial\bar{\bm{x}}}\dd{t}
+\underbrace{\left(
\begin{array}{cc}
    \frac{U^{\text{in}}}{V^*}-\frac{U}{V^*}& 1 \\
     \textit{\textbf{c}}^{\text{in}}-\frac{\textit{\textbf{N}}}{V^*}& 0 \\
  \end{array}\right)}_{{\bar{\bm{g}}}(\bar{\bm{x}})}
  \underbrace{\left(
  \begin{array}{cc}
  q \\
    \dot Q \\
\end{array}\right)}_{{\bar{\bm{u}}}} \dd{t}
+\underbrace{\left(
\begin{array}{c}
0\\
   (\textit{\textbf{z}}-\textit{\textbf{z}}')(\textit{\textbf{r}}_b-\textit{\textbf{r}}_f)\rho_1 V^*\\
  \end{array}\right)}_{\bar{\bm{a}}(\bar{\bm{x}})}\dd{{\omega}_1}
  \\
&&+\underbrace{\left(
\begin{array}{cc}
     \frac{U^{\text{in}}}{V^*}-\frac{U}{V^*}& 1 \\
     \textit{\textbf{c}}^{\text{in}}-\frac{\textit{\textbf{N}}}{V^*}& 0 \\
  \end{array}\right)}_{\bar{\boldsymbol{\gamma}}(\bar{\bm{x}})}
  \underbrace{  \left(
   \begin{array}{cc}
  q  & 0\\
   0& \dot Q \\
\end{array}\right)}_{{\tilde{\bar{\bm{u}}}}}
\underbrace{ \left(
\begin{array}{cc}
\rho_2& 0\\
 0 & \rho_3\\
  \end{array}
\right)}_{{\bar{\boldsymbol{\sigma}}}(\bar{\bm{x}})}
\left(
\begin{array}{c}
\dd{{\omega}_2}\\
\dd{{\omega}_3}\\
  \end{array}
\right)\\
\bar{\bm{y}}&=&\underbrace{\left(
\begin{array}{ccc}
 \frac{U^{\text{in}}}{V^*}-\frac{U}{V^*} & (\textit{\textbf{c}}^{\text{in}}-\frac{\textit{\textbf{N}}}{V^*})^\top \\
 1&0
\end{array}
\right)}_{{\bar{\bm{g}}}^\top(\bar{\bm{x}})}\frac{\partial \bar{H}(\bar{\bm{x}})}{\partial\bar{\bm{x}}}
+
\underbrace{\left(
\begin{array}{cc}
\frac{1}{2}\rho^{2}_2 \frac{M}{\theta}&0\\
0&\frac{1}{2}\rho^{2}_3\frac{1}{\theta}
\end{array}
\right)
}_{\bar{\boldsymbol{\delta}}(\bar{\bm{x}})} \underbrace{\left(
 \begin{array}{c}
 q  \\
    \dot Q \\
\end{array}\right)}_{{\bar{\bm{u}}}}
\end{array}
\right.
\end{split}
\end{equation}
 \hrulefill
 \end{figure*}
\setcounter{equation}{16}
%%%%%%%%%%%%%%%%%%%%%%%%%%%%%%%%%%%%%%%%%%%%%%%%%%%%%%%%%%%%%%%%%%%%%%%%%%%%%%%%%%%%%%%%%%%%%%%%%%%%%%%%%%%%%%%%%%%%%%%%%%%%%%%%%%%%%%%%%%%%%%%%%%%%

\section{Stochastic passivity based control of (state, input)-disturbed CSTR}
The sidCSTR is a sidSPHS if the negative entropy $-S$ is set as the Hamiltonian and a relatively weak condition is added. However, since $-S$ is not bounded from below, the stabilization method derived from the sidSPHS theory \cite{Fang2017} cannot directly apply to the current sidCSTR structure of Eq. (\ref{sid-CSTR}). For this reason, we manage to map the negative entropy to a new Hamiltonian, which is bounded from below and even can behave as a Lyapunov function, and then stabilize the system in probability in this section.

We define the following transformation for the above purpose.
\begin{eqnarray}{\label{ImSGCT}}
\mathscr{T}:~\left\{
\begin{array}{lll}
\bar{\bm{x}}=\textit{\textbf{x}},\\
\bar{H}(\bar{\bm{x}})=-S(\bm{x})+\bm{x}^\top\boldsymbol{\pi}^*\mid_{\textit{\textbf{x}}=\bar{\bm{x}}},\\
\bar{\bm{y}}_d=\textit{\textbf{y}}_d+\boldsymbol{\alpha}(\bm{x})\mid_{\textit{\textbf{x}}=\bar{\bm{x}}},\\
\bar{\bm{u}}=\textit{\textbf{u}},
\end{array}
\right.
\end{eqnarray}
where $\boldsymbol{\pi}^* =\frac{\partial S}{\partial\textit{\textbf{x}}}\mid_{\textit{\textbf{x}}=\bm{x}^*}=\left(\frac{1}{T^*},-\frac{\boldsymbol{\mu}^*}{T^*}^\top\right)^\top$ defined according to Eq. (\ref{EntropyGradient}).

\begin{remark}
The transformed Hamiltonian has explicit physical significance, i.e., representing the availability function \cite{Ydstie1997}. Based on Euler's theorem for homogeneous functions, the entropy function should satisfy $S(\textit{\textbf{x}})=\textit{\textbf{x}}^\top\frac{\partial S}{\partial\textit{\textbf{x}}}=\boldsymbol{\pi}^
\top \textit{\textbf{x}}$. Therefore, the transformed Hamiltonian can be rewritten as
\begin{eqnarray}
\bar{H}(\bar{\bm{x}})=-S(\bar{\bm{x}})+\bar{\boldsymbol{\pi}}^{*\top}(\bar{\bm{x}}-\bar{\bm{x}}^*)+S(\bar{\bm{x}}^*).\notag \end{eqnarray}
Geometrically, it evaluates the distance between the entropy function and its tangent plane at $\bar{\bm{x}}^*$. This together with concavity of the entropy function means $\bar{H}(\bar{\bm{x}})\geq 0$, and moreover, $\bar{H}(\bar{\bm{x}})$ is convex. In the case of setting $V=V^*$, $\bar{H}(\bar{\bm{x}})$ is strictly convex, i.e., its Hessian matrix is positive definite with respect to $\bar{\bm{x}}-\bar{\bm{x}}^{*}$ \cite{Favache2009a}.
\end{remark}

Based on the assumption \textit{A.1}, we fix the volume for the sidCSTR system of Eq. (\ref{sid-CSTR}) to be $V=V^*$, and then apply the transformation of Eq. (\ref{ImSGCT}) to this system.
\begin{proposition}
The transformation defined in Eq. (\ref{ImSGCT}) can transform the sidCSTR system $\Sigma_s$, described by Eq. (\ref{sid-CSTR}), into another sidSPHS in the form of Eq. (\ref{sid-CSTR1}) if we set $\bar{\bm{J}}(\bar{\bm{x}})=\mathbbold{0}_{(p+1)\times (p+1)}$, $\bar{\bm{R}}(\bar{\bm{x}})=\bm{R}(\bm{x}), ~\bar{\bm{g}}(\bar{\bm{x}})=\bm{g}(\bm{x}),~\bar{\bm{u}}=\bm{u}, ~\bar{\bm{a}}(\bar{\bm{x}})=\bm{a}(\bm{x}), ~\bar{\boldsymbol{\gamma}}(\bar{\bm{x}})=\boldsymbol{\gamma}(\bm{x}), ~{\tilde{\bar{\bm{u}}}}=\tilde{\bm{u}}, ~\bar{\boldsymbol{\sigma}}(\bar{\bm{x}})=\boldsymbol{\sigma}(\bm{x})$ and $\bar{\boldsymbol{\delta}}(\bar{\bm{x}})={\boldsymbol{\delta}}(\bm{x})$. Moreover, these two systems are dynamically equivalent.
\end{proposition}

\noindent{\textbf{Proof.}
Since most of variables in the transformed system $\mathscr{T}({\Sigma_s})$ of Eq.(\ref{sid-CSTR1}) are the same as in the original sidCSTR system $\Sigma_s$ of Eq. (\ref{sid-CSTR}), it is clearly that $\mathscr{T}({\Sigma_s})$ is a sidSPHS.

We further prove $\dd \bar{\bm{x}}=\dd \bm{x}$. Since
$$(\bar{\bm{J}}(\bar{\bm{x}})-\bar{\bm{R}}(\bar{\bm{x}}))\frac{\partial \bar{H}(\bar{\bm{x}})}{\partial\bar{\bm{x}}}=-\bm{R}\frac{\partial (-S)}{\partial{\bm{x}}}-\bm{R}\boldsymbol{\pi}^*$$
and
$$\bm{R}\boldsymbol{\pi}^*=\frac{V^*T}{T^*}\sum^{l}_{i=1}\frac{(r_{f_i}-r_{b_i})}{\Delta^\top\textbf{z}_{\cdot i} \boldsymbol{\mu}}
\Delta\textbf{z}_{\cdot i}\Delta^\top\textbf{z}_{\cdot i}\boldsymbol{\mu}^*=\mathbbold{0}_{(p+1)\times 1},$$
where the last equality holds due to $\Delta^\top\textbf{z}_{\cdot i}\boldsymbol{\mu}^*=0$,
we have $$(\bar{\bm{J}}(\bar{\bm{x}})-\bar{\bm{R}}(\bar{\bm{x}}))\frac{\partial \bar{H}(\bar{\bm{x}})}{\partial\bar{\bm{x}}}=(\bm{J}(\bm{x})-\bm{R}(\bm{x}))\frac{\partial (-S)}{\partial\bm{x}},$$ i.e., $\dd \bar{\bm{x}}=\dd \bm{x}$. Therefore, $\mathscr{T}({\Sigma_s})$ and $\Sigma_s$ are dynamically equivalent. $\Box$

\begin{corollary}
Let the sidCSTR systems $\Sigma_s$ and $\mathscr{T}({\Sigma_s})$ be defined by Eqs. (\ref{sid-CSTR}) and (\ref{sid-CSTR1}), respectively, and $\bm{x}^*\in\mathbb{R}^n$ be their equilibrium point. If $\bm{x}^*$ is stable/asymptotically stable in probability for $\mathscr{T}({\Sigma_s})$, then the corresponding result also apply to $\Sigma_s$.
\end{corollary}

This corollary suggests a solution of stabilizing $\Sigma_s$ in probability, i.e., towards stabilizing $\mathscr{T}({\Sigma_s})$ in probability. The latter purpose is easily reached using stochastic passivity theorem \cite{Florchinger1999}.
\begin{theorem}
The transformed sidCSTR system $\mathscr{T}({\Sigma_s})$, governed by Eq. (\ref{sid-CSTR1}), is stochastically passive with respect to $\bar{H}(\bar{\bm{x}})$ if the reaction disturbance $\rho_1$ is small enough such that
\setcounter{equation}{18}
\begin{eqnarray}\label{passivityCond}
 \frac{1}{2}\rho^2_1V^*\sum^p_{j=1}\bm{W}_{j\cdot}\mathbbold{1}_{p}
 \leq
 \sum^l_{i=1}\frac{(r_{f_i}-r_{b_i})(\bm{z}_{\cdot i}-\bm{z}'_{\cdot i})^\top\boldsymbol{\mu}}{T},
 \end{eqnarray}
where
\begin{eqnarray*}
~\left\{
\begin{array}{ll}
\bm{W}=(\frac{\bm{h}\bm{h}^\top}{\theta}-\frac{R}{\textit{\textbf{N}}^\top\mathbbold{1}_{p}}\mathbbold{1}_{p\times p}+\mathrm{diag}\left(\frac{R}{N_j}\right))\circ(\Delta_{\bm{z},\bm{r}}\Delta_{\bm{z},\bm{r}}^\top),\\
\Delta_{\bm{z},\bm{r}}=(\bm{z}-\bm{z}')(\bm{r}_{b}-\bm{r}_{f}).\\
\end{array}
\right.
\end{eqnarray*}
\end{theorem}

\noindent{\textbf{Proof.} Since the transformed system $\mathscr{T}({\Sigma_s})$ in the form of Eq. (\ref{sid-CSTR1}) is a sidSPHS, there exists a necessary and sufficient condition, as \textbf{Theorem 1} says, to suggest this system stochastically passive with respect to its Hamiltonian. Consider
\begin{eqnarray*}
&&\frac{1}{2}\rho^2_1V^{*2}\sum^p_{j=1}\bm
{W}_{j\cdot}\mathbbold{1}_{p}\\
&=&\frac{1}{2}\rho^2_1V^{*2}
\mathrm{tr}\left\{(\frac{\bm{h}\bm{h}^\top}{\theta}-\frac{R\mathbbold{1}_{p\times p}}{\textit{\textbf{N}}^\top\mathbbold{1}_{p}}+\mathrm{diag}(R/N_j))\Delta_{\bm{z},\bm{r}}\Delta_{\bm{z},\bm{r}}^\top\right\}\\
&=&\frac{1}{2}\mathrm{tr}\left\{\frac{\partial^{2} \bar{H}(\bar{\bm{x}})}{\partial {\bar{\bm{x}}^{2}}}\bar{\bm{a}}(\bar{\bm{x}})\bar{\bm{a}}^\top(\bar{\bm{x}})\right\}\\
\end{eqnarray*}
and
\begin{eqnarray*}
\sum^l_{i=1}\frac{(r_{f_i}-r_{b_i})(\bm{z}_{\cdot i}-\bm{z}'_{\cdot i})^\top\boldsymbol{\mu}V^*}{T}=\frac{\partial \bar{H}(\bar{\bm{x}})}{\partial\bar{\bm{x}}}^\top \bar{\bm{R}}(\bar{\bm{x}})\frac{\partial \bar{H}(\bar{\bm{x}})}{\partial \bar{\bm{x}}},
\end{eqnarray*}
which together with Eq. (\ref{passivityCond}) support the condition of Eq. (\ref{passivityCond0})
in \textbf{Theorem 1}. In addition, the current $\bar{\boldsymbol{\delta}}(\bar{\bm{x}})$ defined in Eq. (\ref{sid-CSTR1}) apparently supports the condition of Eq. (\ref{sid-SPHSPassive}) in \textbf{Theorem 1}. Note that $\sum^p_{j=1}\bm
{W}_{j\cdot}\mathbbold{1}_{p}$ has the same sign with $\mathrm{tr}\left\{\frac{\partial^{2} \bar{H}(\bar{\bm{x}})}{\partial {\bar{\bm{x}}^{2}}}\bar{\bm{a}}(\bar{\bm{x}})\bar{\bm{a}}^\top(\bar{\bm{x}})\right\}$, while the latter is nonnegative due to the positive semi-definite $\frac{\partial^{2} \bar{H}(\bar{\bm{x}})}{\partial {\bar{\bm{x}}^{2}}}$ and the positive semi-definite $\bar{\bm{a}}(\bar{\bm{x}})\bar{\bm{a}}^\top(\bar{\bm{x}})$. Therefore, as long as the reaction disturbance $\rho_1$ is small enough, the transformed sidCSTR system can be stochastically passive with respect to $\bar{H}(\bar{\bm{x}})$. $\Box$

%%%%%%%%%%%%%%%%%%%%%%%%%%%%%%%%%%%%%%%%%%%%%%%%%%%%%%%%%%%%%%%%%%%%%%%%%%%%%%%%%%%%%%%%%%%%%%%%%%%%%%%%%%%%%%%%%%%%%%%%%%%%%%%%%%%%%%%%%%%%%%%
\newcounter{TempEqCnt2}                         % 创建临时变量TempEqCnt
\setcounter{equation}{20}% 当前公式序号变为x，x等于长公式应有的序号减1.
\begin{figure*}[!t]
\hrulefill
\begin{equation}{\label{S_1}}
\begin{split}
~\left\{
\begin{array}{lll}
\dd{\textit{\textbf{x}}}&= &
-\underbrace{\frac{(k_{0f}N_A e^{-\frac{E_{f}}{RT}}-k_{0b}N_B e^{-\frac{E_{b}}{RT}})T}{\mu_A-\mu_B}\left(
 \begin{array}{ccc}
    0 &0& 0\\
  0& 1& -1\\
      0& -1& 1\\
\end{array}
\right)}_{{\textit{\textbf{R}}}(\textit{\textbf{x}})}
\frac{\partial (-S)}{\partial\textit{\textbf{x}}} \dd{t}
+\underbrace{\left(
\begin{array}{cc}
c^{\mathrm{in}}_A {h}^{\mathrm{in}}-(\frac{N_A}{V}{h}_A+\frac{N_B}{V}{h}_B)& 1 \\
    c^{\mathrm{in}}_A-\frac{N_A}{V}& 0 \\
    \frac{-N_B}{V}& 0\\
  \end{array}\right)}_{{\textit{\textbf{g}}}(\textit{\textbf{x}})}
  \underbrace{\left(
  \begin{array}{cc}
  q \\
    \dot Q \\
\end{array}\right)}_{{\textit{\textbf{u}}}} \dd{t}
\\
&&+\underbrace{\left(
\begin{array}{c}
   0\\
   -\rho_1(k_{0f}N_A e^{-\frac{E_{f}}{RT}}-k_{0b}N_B e^{-\frac{E_{b}}{RT}})\\
   \rho_1(k_{0f}N_A e^{-\frac{E_{f}}{RT}}-k_{0b}N_B e^{-\frac{E_{b}}{RT}})\\
  \end{array}\right)}_{{\bm{a}(\bm{x})}}\dd{{\omega}_1}
+\underbrace{\left(
\begin{array}{cc}
c^{\mathrm{in}}_A{h}^{\mathrm{in}}-(\frac{N_A}{V}{h}_A+\frac{N_B}{V}{h}_B)& 1 \\
c^{\mathrm{in}}_A-\frac{N_A}{V}& 0 \\
    \frac{-N_B}{V}& 0\\
  \end{array}\right)}_{{\boldsymbol{\gamma}(\bm{x})}}
  \underbrace{  \left(
   \begin{array}{cc}
  q  & 0\\
   0& \dot Q \\
\end{array}\right)}_{{\tilde{\textit{\textbf{u}}}}}
\underbrace{ \left(
\begin{array}{cc}
{\rho}_2 & 0\\
 0 & {\rho}_3\\
  \end{array}
\right)}_{{\boldsymbol{\sigma}}(\bm{x})}
\left(
\begin{array}{c}
\dd{{\omega}_2}\\
\dd{{\omega}_3}\\
  \end{array}
\right)\\
\textit{\textbf{y}}&=&\underbrace{\left(
\begin{array}{ccc}
  c^{\mathrm{in}}_A{h}^{\mathrm{in}}-(\frac{N_A}{V}{h}_A+\frac{N_B}{V}{h}_B)&c^{\mathrm{in}}_A-\frac{N_A}{V}&\frac{-N_B}{V}  \\
   1&   0&  0\\
  \end{array}
\right)}_{{\textit{\textbf{g}}}^\top(\textit{\textbf{x}})}\frac{\partial  (-S)}{\partial\textit{\textbf{x}}}
+
\underbrace{\left(
\begin{array}{cc}
\frac{1}{2}\rho^{2}_2 \frac{M}{\theta}&0\\
0&\frac{1}{2}\rho^{2}_3\frac{1}{\theta}\\
\end{array}
\right)}_
{\boldsymbol{\delta}(\bm{x})}
\underbrace{\left(
 \begin{array}{c}
 q  \\
    \dot Q \\
\end{array}\right)}_{{\textit{\textbf{u}}}}
\end{array}
\right.
\end{split}
\end{equation}
\hrulefill
 \end{figure*}
\setcounter{equation}{19}
%%%%%%%%%%%%%%%%%%%%%%%%%%%%%%%%%%%%%%%%%%%%%%%%%%%%%%%%%%%%%%%%%%%%%%%%%%%%%%%%%%%%%%%%%%%%%%%%%%%%%%%%%%%%%%%%%%%%%%%%%%%%%%%%%%%%%%%%%%%%%%%

\begin{remark}
The condition of Eq. (\ref{passivityCond}) is not so harsh from the viewpoint of practice since the order of magnitude of $\rho_1^2$ is quite low compared with other macroscopic variables, which leads to Eq. (\ref{passivityCond}) to be most likely true in the practical cases.
\end{remark}

The stochastic passivity of $\mathscr{T}({\Sigma_s})$ can render a simple and effective way to stabilize it in probability, and further for the original sidCSTR system $\Sigma_s$.

\begin{theorem}
For any sidCSTR system $\Sigma_s$ given by Eq. (\ref{sid-CSTR}), let $\bm{x}^*$ be one of its equilibrium points. Assume that its transformed version $\mathscr{T}(\Sigma_s)$ modeled by Eq. (\ref{sid-CSTR1}) is stochastically passive with respect to $\bar{H}(\bar{\bm{x}})$, and there exists a positive definite matrix $\textbf{K}\in \mathbb{R}^{2\times 2}$ such that $\textit{\textbf{I}}+\textbf{K}\bar{\boldsymbol{\delta}}$ is invertible.  Then the controller in the form of
\begin{eqnarray}\label{Controller}
\bar{\bm{u}}=-\big[\bm{I}+\textbf{K}\bar{\boldsymbol{\delta}}\big]^{-1}\textbf{K}
\bar{\bm{g}}^\top(\bar{\bm{x}})\frac{\partial
\bar{H}(\bar{\bm{x}})}{\partial {\bar{\bm{x}}}}
\end{eqnarray}
can locally asymptotically stabilize $\Sigma_s$ at $\bm{x}^*$ in probability if they are connected in negative feedback.
\end{theorem}

\noindent{\textbf{Proof.}
Since the transformed system $\mathscr{T}(\Sigma_s)$ modeled by Eq. (\ref{sid-CSTR1}) is stochastically passive with respect to $\bar{H}(\bar{\bm{x}})$, we have $\frac{\dd \bar{H}(\bar{\bm{x}})}{\dd t}\leq \bar {\bm{y}}^\top \bar{\bm{u}}$. Combining Eqs. (\ref{sid-CSTR1}) and (\ref{Controller}) yields
$$\bar{\bm{u}}=-\boldsymbol{K}(\bar{\bm{g}}^\top(\bar{\bm{x}})\frac{\partial\bar{H}(\bar{\bm{x}})}{\partial\bar{\bm{x}}}+\boldsymbol{\bar{\delta}} {\bar{\bm{u}}})=-\boldsymbol{K}\bar {\bm{y}}.$$
Hence, we obtain $\frac{\dd \bar{H}(\bar{\bm{x}})}{\dd t}\leq \bar {\bm{y}}^\top \bar{\bm{u}}=-\bar {\bm{y}}^\top \boldsymbol{K}\bar {\bm{y}}\leq 0$ with equality hold if and only if $\bar{\bm{y}}=\mathbbold{0}_2$, i.e., $\bar{\bm{x}}=\bar{\bm{x}}^*$. Since $\bar{H}(\bar{\bm{x}})$ is positive definite with respect to $\bar{\bm{x}}-\bar{\bm{x}}^*$, the controller can locally asymptotically stabilize $\mathscr{T}(\Sigma_s)$ at $\bar{\bm{x}}^*$ in probability. This together with \textbf{Corollary 1} suggests that the equilibrium point $\bm{x}^*$ for $\Sigma_s$ can be locally asymptotically stabilized in probability by the proposed controller. $\Box$

%%%%%%%%%%%%%%%%%%%%%%%%%%%%%%%%%%%%%%%%%%%%%%%%%%%%%%%%%%%%%%%%
\section{Case study}
In this section, the proposed port-Hamiltonian method is applied to stabilizing in probability a CSTR system \cite{Hoang2011} with state and input both disturbed, in which a first-order reversible chemical reaction $A\rightleftarrows B$ takes place, and only the pure $A$ component is fed into the reactor. The system follows all assumptions $A.1-A.5$.

The dynamical equations of this system ${\Sigma_s}$ may be written in the sense of It$\mathrm{\hat{o}}$ to be
\begin{eqnarray*}
\left\{
\begin{array}{lll}
\dd{U}&=&q[c^{\mathrm{in}}_A{h}^{\mathrm{in}}-(\frac{N_A}{V}{h}_A+\frac{N_B}{V}{h}_B)](\dd{t}+{\rho_2}\dd{\omega_2})\\
&&+\dot Q(\dd{t}+{\rho_3}\dd{\omega_3}),
\\
\dd{N_A}&=&(k_{0b}N_B e^{-\frac{E_{b}}{RT}}-k_{0f}N_A e^{-\frac{E_f}{RT}}) (\dd{t}+{\rho_1}\dd{\omega_1})
\\
&&q( c^{\mathrm{in}}_A-\frac{N_A}{V})(\dd{t}+{\rho_2}\dd{\omega_2}),
\\
\dd{N_B}&= &(k_{0f}N_A e^{-\frac{E_{f}}{RT}}-k_{0b}N_B e^{-\frac{E_{b}}{RT}})(\dd{t}+{\rho_1}\dd{\omega_1})
\\
&&-q\frac{N_B}{V}(\dd{t}+ {\rho_2}\dd{\omega_2}),
\end{array}
\right.
\end{eqnarray*}
where some related parameters values are listed in Table $1$. For this system, there are up to three equilibria with the middle one unstable while other two stable. The current task focuses on stabilizing in probability the middle unstable equilibrium point, which is denoted by $\bm{x}^*=(U^*,N_A^*,N_B^*)^\top$. The reactor works under the initial conditions given in Table $2$, based on which it is easy to calculate the concerning equilibrium point to be $(1157.5,1.3,0.7)^\top$, and the corresponding $V^*=0.001$m$^3$ and $T^*=331.9$K. We also exhibit them in Table $2$.

\begin{table}[hbp]
\tiny
\begin{tabular}{p{0.6cm}p{1.9cm}p{4.5cm}}
\multicolumn{3}{c}{\small{Table 1. Parameters of the sidCSTR system: $A\rightleftarrows B$}} \\\hline
Symbol&Numerical value&Variable name\\\hline
$C_{P_A}$&75.24J/K/mol& Heat capacity of species $A$\\
$C_{P_B}$&60J/K/mol& Heat capacity of species $B$\\
$h_{A_{\mathrm{ref}}}$&0J/mol&Reference enthalpy of $A$\\
$h_{B_{\mathrm{ref}}}$&-4575J/mol&Reference enthalpy of $B$\\
$k_{0f}$&$0.12\times10^{10}$/s& Forward kinetic constant\\
$E_{f}$&72331.8J/mol&Forward reaction activation energy\\
$k_{0b}$&$1.33\times10^8$/s& Backward kinetic constant\\
$E_{b}$&74826J/mol& Backward reaction activation energy\\
$P$&$10^5$Pa&Pressure\\
$R$&$8.314$ J/K/mol &Molar gas constant\\
$T_{\mathrm{ref}}$&$300$K&Reference Temperture\\
$s_{A_{\mathrm{ref}}}$&$50.6$J/K/mol &Reference entropy of $A$\\
$s_{B_{\mathrm{ref}}}$&$180.2$J/K/mol &Reference entropy of $B$\\
$\rho_1,\rho_3$&$0.1, 0.05$&Reaction and heat exchange disturbance\\
$\rho_2$ & $5\times 10^{-7}$ & Inlet/outlet flow disturbance \\
$\lambda$&$0.05808$J/K/s&Heat transfer coefficient\\
\hline
\end{tabular}
\end{table}

\begin{table}[hbp]
\tiny
\begin{tabular}{p{1.8cm}p{2.0cm}p{1cm}p{1.9cm}}
\multicolumn{4}{c}{\small{Table 2. Initial conditions and setpoint values}}\\
\hline
Initial&Value&Setpoint& Value\\\hline
$T^{\text{in}},T(0),T_w(0)$ & $310,342,299.97$K  & $U^*,T^*$ & $1157.5$J,$331.9$K\\
  $q(0)$ & $9.15\times 10^{-4}$m$^3$/s & $V^*$ & $0.001$m$^3$\\
  $N_{\text{A}}(0),N_{\text{B}}(0)$ & $1,1$mol & $N_A^*$, $N_B^*$ &  $1.3$, $0.7$mol\\
 $c^{\text{in}}_A$  & $2000$mol/m$^3$/s & $q^*$ &$9.15\times 10^{-6}$m$^3$/s\\
\hline
\end{tabular}
\end{table}

By setting $H(\bm{x})=-S$ and then applying \textbf{Proposition 2}, we rewrite the above dynamical equation as the form of sidSPHS, shown in Eq. (\ref{S_1}), where
\setcounter{equation}{21}
\begin{eqnarray}{\label{delta2}}
~\left\{
\begin{array}{lll}
\frac{M}{\theta}=\frac{{c^{\mathrm{in}}_A}^2 \textit {R}}{(N_A+N_B)N_{A} N_{B}}+\frac{(c^{\mathrm{in}}_A{h}_A-c^{\mathrm{in}}_A{h}^{\mathrm{in}})^2}{\theta},\\
\theta=(N_A C_{P_A}+N_B C_{P_B})T^2. \\
 \end{array}
\right.
\end{eqnarray}
It is clear that at given $\rho_2,\rho_3$ Eq. (\ref{norm cond}) is true to support Eq. (\ref{S_1}) to be a sidSPHS. Further applying \textbf{Proposition 3} will produce a transformed but equivalent sidSPHS to the original one of Eq. (\ref{S_1}). Moreover, at the given $\rho_1$, the condition of Eq. (\ref{passivityCond}) in \textbf{Theorem 2} holds, which means the transformed sidCSTR system stochastically passive with respect to its Hamiltonian. Finally, based on \textbf{Theorem 3}, we set the proportional gain matrix $\boldsymbol{K}$ to be
$\boldsymbol{K}=\text{diag}(K_1, K_2)=\text{diag}(1.64\times 10^{-7}, 27430)$, which supports $\textit{\textbf{I}}+\textbf{K}\bar{\boldsymbol{\delta}}$ to be invertible, and then design the following controller utilizing Eq. (\ref{Controller})
\begin{eqnarray*}{\label{controller1}}
~\left\{
 \begin{array}{lll}
 q&=&-K_1 (1+ \frac{1}{2}{\rho_2}^2\frac{M}{\theta}K_1)^{-1}[- \frac{N_B}{V}(-\frac{\mu^*_B}{T^*}+\frac{\mu_B}{T})\\
 &&+(c^{\mathrm{in}}_A{h}^{\mathrm{in}}-\frac{N_A}{V}{h}_A-\frac{N_B}{V}{h}_B )(\frac{1}{T^*}-\frac{1}{T})\\
&&+(c^{\mathrm{in}}_A-\frac{N_A}{V})(-\frac{\mu^*_A}{T^*}+\frac{\mu_A}{T})],\\
 \dot Q&=&-K_2 (1+ \frac{1}{2}{\rho_3}^2\frac{1}{\theta}K_2)^{-1} (\frac{1}{T^*}-\frac{1}{T})
 \end{array}
 \right.
 \end{eqnarray*}
to locally asymptotically stabilize the considered sidCSTR system in probability. The manipulated variable $T_w$ may be further calculated by combining the above second equation with the expression $\dot Q=\lambda (T_w-T)$ as
\begin{equation*}
 T_w=\frac{-K_2 (1+ \frac{1}{2}{\rho_3}^2\frac{1}{\theta}K_2)^{-1} (\frac{1}{T^*}-\frac{1}{T})}{\lambda}+T.
\end{equation*}

To observe the state converging behaviors for this system, Figs. $2$ and $3$ exhibit the response of the temperature and molar mass to time, respectively. As one can expect, the state vector $(T,N_A,N_B)^\top$ can asymptotically converge to the equilibrium state $(T^*,N_A^*,N_B^*)^\top$ as the controller is put into force. Basically, from $t=3$s the converging behaviors happen. These phenomena suggest that the proposed stochastic passivity based controller is qualified for stabilizing the sidCSTR system in probability. We also present the evolution of the manipulated variables $q$ and $T_w$ in Figs. $4$ and $5$. In like manner, they will asymptotically approaching the stable values $q\approx 0$ and $T_w\approx T^*$ as the state tends to the equilibrium point.

\begin{figure}
\begin{center}
{\includegraphics[width=3.3in]{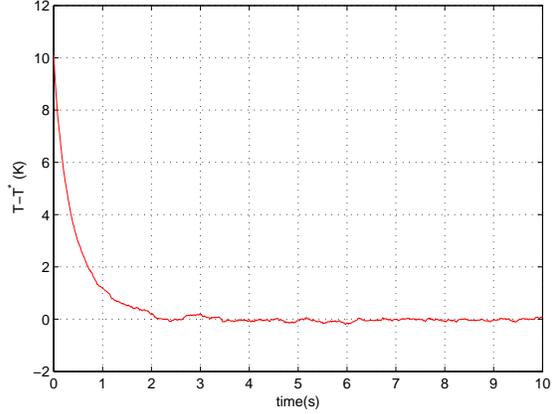}} {\caption{
Temperature response of the sidCSTR system.}}
\end{center}
\end{figure}

\begin{figure}
\begin{center}
{\includegraphics[width=3.3in]{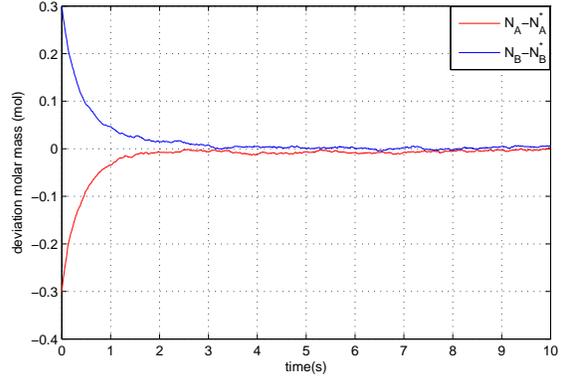}} {\caption{Molar mass responses of the sidCSTR system.}}
\end{center}
\end{figure}

\begin{figure}
\begin{center}
{\includegraphics[width=3.3in]{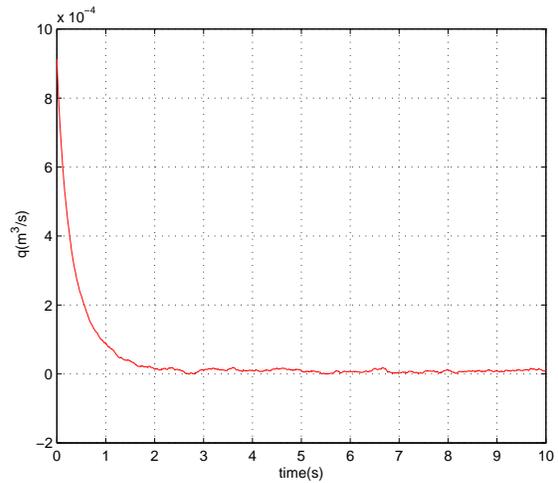}} {\caption{Input/output volume flow evolution.}}
\end{center}
\end{figure}

\begin{figure}
\begin{center}
{\includegraphics[width=3.3in]{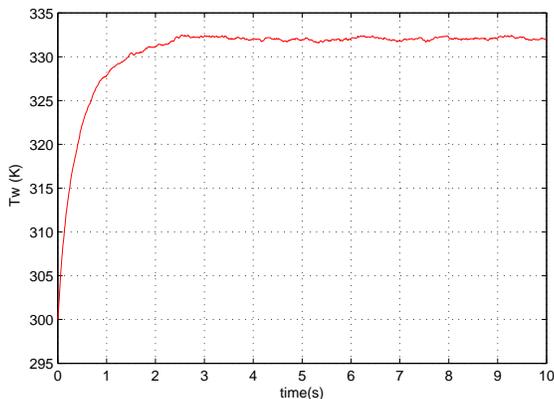}} {\caption{Jacket temperature evolution.}}
\end{center}
\end{figure}

%%%%%%%%%%%%%%%%%%%%%%%%%%%%%%%%%%%%%%%%%%%%%%%%%%%%%%%%%%%%%%%%
\section{Conclusion}
In this paper, a port-Hamiltonian based method is proposed to model and stabilize in probability the CSTR systems with state and input both disturbed. The following conclusions are reached:

(i) The sidCSTR systems may be formulated into a sidSPHS structure if the opposite entropy function is taken as the Hamiltonian. Moreover, the degenerated deterministic version can represent the first law and the second law of thermodynamics simultaneously.

(ii) The sidCSTR systems have the property of stochastic passivity with respect to the Hamiltonian when the Hamiltonian is transformed from the opposite entropy function to the availability function.

(iii) The sidCSTR systems can be locally asymptotically stabilized if a stochastic passivity based controller is connected in negative feedback.
%%%%%%%%%%%%%%%%%%%%%%%%%%%%%%%%%%%%%%%%%%%%%%%%%%%%%%%%%%%%%%%%

\begin{ack}                               % Place acknowledgements
This work was supported by
the National Natural Science Foundation of China under Grant No. 11671418,
11271326 and 61611130124, and the Research Fund for the Doctoral Program
of Higher Education of China under Grant No. 20130101110040.
\end{ack}

 %\bibliographystyle{plain}
 %\bibliography{mybib}

\end{document}